\documentclass[a4paper]{article}
\usepackage{geometry}
\usepackage{ulem}
\usepackage{graphicx}
\usepackage{amsmath, amsfonts}
\usepackage{titlesec}
\usepackage{booktabs}
\usepackage{placeins}
\usepackage{caption}
\title{\textbf{Improved Fast Iterative Algorithm for Eikonal Equation for GPU Computing}}
\author{Yuhao Huang $^a$\\\small$^a$ Department of Engineering Science and Applied Mathematics, \\\small Northwestern University, Evanston, 60201, USA}
\begin{document}
\begin{titlepage}
\maketitle 
\begin{center}
\textbf{Abstract}
\end{center}
\begin{center}
\begin{minipage}{0.8\textwidth}

In this paper we propose an improved fast iterative method to solve the Eikonal equation, which can be implemented in parallel. We improve the fast iterative method for Eikonal equation in two novel ways, in the value update and in the error correction.  The new value update is very similar to the fast iterative method in that we selectively update the points, chosen by a convergence measure, in the active list. However, in order to reduce running time, the improved algorithm does not run a convergence check of the neighboring points of the narrow band as the fast iterative method usually does. The additional error correction step is to correct the errors that the previous value update step may cause. The error correction step consists of finding and recalculating the point values in a separate remedy list which is quite easy to implement on a GPU. In contrast to the fast marching method and the fast sweeping method for the Eikonal equation, our improved method does not need to compute the solution with any special ordering in neither the remedy list nor the active list.  Therefore, our algorithm can be implemented in parallel. In our experiments, we implemente our new algorithm in parallel on a GPU and compare the elapsed time with other current algorithms. The improved fast iterative method runs faster than the other algorithms in most cases through our numercal studies.\\\\\\

\end{minipage}
\end{center}

\section{Introduction}
Moving interfaces occur in a wide range of scientific and engineering applications. Some of them are obvious such as bubbles rising in a fluid, crystals solidifying in an aqueous solution while some of them are not quite obvious like optimizing solutions for robotic path planning. The numerical simulation methods in literature have been applied in  the fields of computer vision, image processing, geoscience, medical imaging and robotic path planning (see more details in [9, 10, 11, 12]). 
\\\\
Let $\Omega$ be a spatial domain in which an interface $\Gamma$ of co-dimension one will be contained. Let $\textbf n$ be the unit normal to $\Gamma$ at \textbf{x}. For a moving interface, for each $\textbf x \in \Gamma$ on the interface, there is a prescribed velocity vector $\textbf v$ which is specific to the application and might be determined from many different places.  For example, it may depend on the local geometry, prescribed velocity potential, or other formulae derived from the physics of the application. The velocity vector $\textbf v$ can be broken down into its normal and tangential components $\textbf v = \textbf v^{\perp} + \textbf v^{\parallel} $ respectively where $\textbf v^\perp=(\textbf v\cdot \textbf n)\textbf n$ is the normal comoponent and $\textbf v^\parallel=\textbf v-\textbf v^\perp$ is the component parallel to $\Gamma$. Figure~\ref{sketch} illustrates these vectors relative to the interface.\\\\
\begin{figure}[h]
$$\includegraphics[scale=0.4]{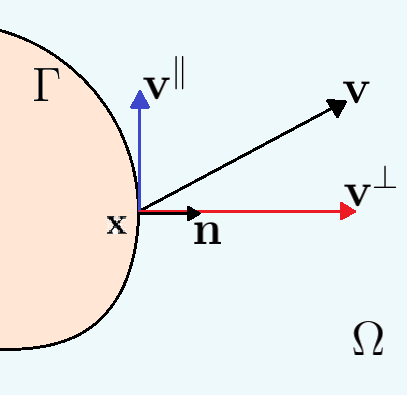}$$
\caption{An illustration of the interface and relative vectors.}
\label{sketch}
\end{figure}
\noindent If $\textbf v^{\perp} = 0$ for each $\textbf x \in \Gamma$, then the interface will remain stationary because  the entire flow is tangential to the interface. The only part of the velocity which contributes to the motion of the interface is $\textbf v^{\perp}$. We will call $\textbf v^{\perp}$ the flow in the normal direction and let $\textbf v^{\perp} = F\textbf n$ for some speed function $F($\textbf x$)$. We define the basic flow equation for interface motion as, 
$$\frac{d\textbf x}{dt} = F(\textbf x, t)\textbf n,\eqno(1)$$
where in this paper we will assume $F(\mathbf x,t)>0$. With this assumption, the interface can be represented implicitly by a function $\phi(\mathbf x)$, where the set $\phi(\mathbf x)=t$ is the location of the interface at time t. The function $\phi$ is often referred to as the time of crossing map. To construct $\phi$, suppose $\textbf x \in \Gamma$, then at time $t$,
$$t=\phi(\textbf x(t)).\eqno(2)$$
Differentiating with respect to t gives,
$$1=\nabla \phi  \cdot  \textbf x'(t) = \nabla \phi \cdot (F\textbf n).$$
Since $t = \phi(\textbf x)$ is the location of the interface at time $t$, and this is an isocontour, then it follows that,
$$\textbf n = \frac{\nabla \phi}{|| \nabla \phi ||}.$$
We thus obtain,
$$1 = \nabla \phi \cdot ( F\frac{\nabla \phi}{|| \nabla \phi ||} ) = F||\nabla \phi||.$$
$$\Longrightarrow ||\nabla \phi|| = \frac{1}{F},\eqno(3)$$
which is the Eikonal equation, a special case of nonlinear Hamilton-Jacobi partial differential equations.  \\\\
To solve this equation, the Godunov upwind finite difference scheme is widely used. For example, the first order Godunov upwind discretization of (3) on a 2D
uniform grid is given by,
$$\left[\frac{(\phi_{i,j}-\phi_{i,j}^{xmin})^+}{\Delta x}\right]^2+\left[\frac{(\phi_{i,j}-\phi_{i,j}^{ymin})^+}{\Delta y}\right]^2=\frac{1}{F_{i,j}^2},\eqno(4) $$
where $\phi_{i,j}$ is the discrete approximation to $\phi$ at $\textbf x = (x_i,y_j)$, and $\phi_{i,j}^{xmin}=\min (\phi_{i-1,j},\phi_{i+1,j})$, $\phi_{i,j}^{ymin}=\min (\phi_{i,j-1},\phi_{i,j+1})$, and $(x)^+=\max(x,0)$. This can be a good local solver for the Eikonal equation.\\\\
There are many numerical strategies for solving (1). One method is to consider this problem as a stationary boundary value problem and then to design an efficient numerical algorithm to solve the system of nonlinear equations. For example, the fast marching method [1, 2, 3] is of this type. In the fast marching method, the update of the solution follows the causality that the solution is updated  grid point by grid point in order of increasing values of $\phi$. To accomplish this efficiently, a heapsort queue and an upwind finite difference discretization scheme such as (4) is used. The complexity for the fast marching method is $O(N\log N)$ for $N$ gridpoints, where the $\log N$ factor is due to maintaining the heapsort. The fast sweeping method [4, 5] is also a commonly used method to solve the system based on the same Godunov upwind discretization scheme and using Gauss-Seidel iterations with alternating sweep ordering. The fast sweeping method follows the causality along characteristics in a parallel way; i.e., all characteristics are divided into a finite number of groups according to their directions and each Gauss-Seidel iteration with a specific sweep direction covers a group of characteristics simultaneously. The fast sweeping method is extremely simple to implement and the complexity of the algorithm is $O(N)$ for a total of $N$ grid points [4].\\\\
Parallel implementation of the fast marching and fast sweeping algorithms is challenging due to their inherent structure. For the fast marching algorithm, the heapsort does not easily allow for massively parallel solutions, such as those available with single instruction multiple datastream (SIMD) architectures. Furthermore, grid points must be updated one at a time so that a parallel implementation is not readily apparent. The fast sweeping method requires reading from and writing to the single memory location that stores the grid value which is inefficient or prohibited on some of the most efficient parallel architectures. The fast iterative method (FIM) [6] was developed to solve the problem on a GPU which meets the requirement of SIMD. However, in the FIM, each iteration, the neighboring points of the narrow band points require a convergence check by calculating the solution for Equation (2), which is time-consuming even on a GPU.\\\\
In this paper, we present and analyze an improved iterative algorithm based on the fast iterative method. There are novel modifications to the current algorithm, one for updating values and one for error correction. The new value update part is very similar to the fast iterative method in that we selectively update the points in the active list chosen by a convergence measure but do not check the convergence of the neighboring points of the narrow band in each iteration. The error correction is to get rid of the errors that the previous value update step may cause. The error correction is to find and recalculate the point values in a separate remedy list which is quite easy to implement on a GPU. The resulting complexity of the improved fast iterative method is the same as the fast iterative method $\textbf{---}$ $O(N)$ [6]. We can update all the point values in the active list and then remedy the point values in the remedy list in parallel on a GPU, which saves significant time.\\\\
This paper proceeds as follows. In the next section, we introduce previous work such as the fast marching method, the fast sweeping method and the fast iterative method.  In Section $3$, we introduce the proposed improved fast iterative method in general. In Section $4$, we show the implementation in detail using CUDA and numerical results on a number of different examples and compare with the existing methods. In section $5$, we summarize the paper and discuss the future research directions related to this work.

\section{Previous Work}
In this section, we discuss three numerical methods to solve the Eikonal equation, namely the fast marching method, the fast sweeping method and the fast iterative method. Each of them represent different strategies and logic. In this discussion, we will introduce these algorithms in detail as well as their relevant background. \\

\subsection{The Fast Marching Method}	
The fast marching method (FMM) is a novel method for capturing the motion of interfaces. It was introduced by Sethian [1]. Where applicable, it is vastly faster than any other numerical method by an order of magnitude. It is also an interesting method in its own right as well as a key piece used in modern level set method implementations. It uses a heap to sort points on the moving wavefront, which is based on the property of the solution to guarantee the characteristics are a steepest descent on $\phi$. In this algorithm, each point value only depends on its neighboring points whose values are smaller and the point value is updated in the order maintained in a heap-sort queue. Here, two points $\mathbf{x}_{i,j}$, $\mathbf{x}_{k,\ell}$ are neighbors if $|i-k|+|j-\ell|=1$.\\\\
There are three different regions defined by the fast marching method, the $Accepted$, $Tentative$ and $Distant$ regions. The point values in the $Accepted$ region are already found and fixed. The points that neighbor points in the $Accepted$ region, but are not themselves in the $Accepted$ region are considered to be in the $Tentative$ region. The remaining points outside of the boundary are considered in the $Distant$ region and their values are infinite and not updated until added into the $Tentative$ region. The algorithm steps are as follows:\\

$\bullet$ At the beginning, the algorithm puts the boundary conditions into the $Accepted$ set, then initializes the \textit{Tentative} set to be all the points that are neighbors of points in the \textit{Accepted} set but not in the \textit{Accepted} set. For each point in the \textit{Tentative} set, the value is set by solving Equation (4). The \textit{Tentative} set is constructed as a min-heap structure so that the top of the heap always has the point with the lowest $\phi$ value. All the remaining points are put in the \textit{Distant} set.\\

$\bullet$ Each iteration, the point on the top of the heap queue is removed from the \textit{Tentative} set and moved into the $Accepted$ set. The values of its neighbor points that are not in the \textit{Accepted} set are updated by solving Equation (4) and the neighbor points that were in the \textit{Distant} set are moved into the $Tentative$ set. This find-and-update procedure continues until the $Tentative$ set is empty. \\\\
The computational cost of the fast marching method is $O(N \log N)$ where $N$ is the total number of grid points. The log $N$ part arises due to the cost of sorting the \textit{Tentative} set.\\

\subsection{The Fast Sweeping Method}
Zhao [4] proposed the fast sweeping method (FSM) as an alternative to the fast marching method. It is an iterative method using upwind differencing for discretization and using Gauss-Seidel iterations with alternating sweep ordering to solve the discretized system based on the idea that each sweep follows a family of characteristics of the corresponding Eikonal equation in a certain direction simultaneously. The Eikonal equation is solved on an n-dimensional grid using at least $2^n$ directional sweeps, one per each quadrant, within a Gauss-Seidel update scheme.\\\\
The fast sweeping method follows the Godunov upwind difference scheme [7] to discretize the partial differential equation at the points on the boundary so point values only depend on neighboring points whose values are smaller. The algorithm steps are as follows:\\

$\bullet$ For initialization, grid points at the boundary are assigned values according to the boundary conditions and large positive values are assigned at all other grid points. \\

$\bullet$ Each iteration, for the points whose values are not fixed during the initialization, it computes the solutions based on the upwind difference scheme by solving Equation (4) with alternating sweeping orderings repeatedly:
$$1)\quad i=1, 2, \ldots,I,\quad j=1, 2, \ldots, J,\qquad \qquad \qquad2)\quad i=I, I-1, \ldots, 1,\quad j=1,2,\ldots,J,$$
$$3)\quad i=I,I-1,\ldots,1,\quad j=J,J-1,\ldots,1,\qquad \quad 4)\quad i=1,2,\ldots,I,\quad j=J,J-1,\ldots,1.$$
The sweeping repeats until convergence. \\\\
The algorithm is optimal in the sense that a finite number of iterations is needed. So the complexity of the algorithm is $O(N)$ for a total of $N$ grid points. \\

\subsection{The Fast Iterative Method}
The fast iterative method (FIM) proposed by Won-Ki Jeong and Ross Whitaker [6, 8] is an Eikonal solver based on a selective iterative method. By using this method, it is much easier to produce good overall performance, cache coherence, and scalability across multiple processors. Compared with the former two algorithms, the advantages of the fast iterative method are: \\
\hspace*{2em}$\bullet$ it does not require a particular update order\\
\hspace*{2em}$\bullet$ It does not require a separate, sorted, heterogeneous data structure\\
\hspace*{2em}$\bullet$ It's possible to update multiple points simultaneously\\\\
In the fast iterative method,  it maintains a narrow band, called the $Active$ list, for storing the grid points that are being updated. The main algorithm consists of two parts, the initialization step and the update step. \\

$\bullet$ For initialization,  it sets the initial values for grid points located in the areas specified by the boundary conditions and labels them as source points. The remaining point values are set to infinity. Then, the neighbors of the source points are added to the $Active$ list. \\

$\bullet$ In the update step, the values of points in the $Active$ list are updated depending on their neighboring source points. Once the value of an active point converges, it is added to the source point set and then there is a check whether its neighboring point values have converged or not. If a neighboring point is infinity or is not converged, it is added to the $Active$ list. The iteration stops when there are no points in the $Active$ list. \\\\
This method meets the requirement of SIMD so that it can be implemented on a GPU to solve the problem. We can update the $\phi$ values of the points in the $Active$ list simultaneously! \\\\
However, the fast iterative method needs to check whether the points of the $Active$ list are convergent or not by solving Equation (4). This is because the values of the neighboring points in the $Active$ list may change after the values in the $Active$ list are updated, which means convergence checking of the active points can't be done until after the update step for the $Active$ list is completed in each iteration. This convergence checking step can be time-consuming. Even in a parallel implementation, each update of the value of a grid point requires a check whether the neighboring point values have changed or not. These two steps cannot be implemented in parallel because of the causal relationship. In most cases, each point in the $Active$ list has more than two neighboring points that are not in the $Active$ list and their values need to be calculated in each iteration so this can be inefficient.

\section{Improved Fast Iterative Method}
In this section, we propose a new method based on the fast iterative method which also meets the requirement for SIMD. It also produces good overall performance, cache coherence, and scalability across multiple processors. Similar to FIM, the new algorithm does not impose a special update order so it does not use any data structure for sorting points. A special update order requires random access of memory that will destroy cache coherency. In the new algorithm, multiple points can be updated simultaneously enabling the algorithm to fully utilize a GPU parallel architecture.\\

\subsection{Introduction to the Improved FIM}
The improved fast iterative method is a modification of the fast iterative method [6]. It changes some iteration rules to improve the speed and adds a remedy algorithm to correct the possible errors that may happen in the update step. The fast iterative method [6] was based on observations from two numerical solvers. One is an iterative method proposed by Rouy etal. [8] that updates all point values in each iteration. This algorithm is simple to implement but it does not rely on a causality principle and hence every grid point must be visited in each iteration making it inefficient. The other is the fast marching method which updates the point values selectively from a heap queue but does not support parallel computing. \\\\
As opposed to the fast marching method, the improved fast iterative method updates the point values without maintaining expensive data structures. It uses three sets for the points being updated. The first set is called the $Active$ set where the grid points being updated are stored. In each iteration, convergent points are removed and new points are added into this set. The second set is called the $Source$ set that collects the points whose values have been updated in the update step. After the $Active$ set becomes empty, the update step of the improved FIM ends and the remedy step begins where it maintains an additional set called the $Remedy$ set. When the update step is completed, the improved FIM updates all the grid point values once and then adds the points whose values change into the $Remedy$ set. In each iteration, it updates the values of points in the $Remedy$ set and points adjacent to this set. In each iteration of the remedy step, if a point not in the Remedy set has its value change, then it is added to the $Remedy$ set. If the value of a point in the Remedy set doesn't change, then it is removed from the Remedy set. The remedy step continues until the $Remedy$ set is empty. Although it is very expensive to do the remedy work serially, it is easy to parallelize to run on a GPU.\\

\subsubsection{Local Solver}
We use the Godunov upwind discretization to update the values of the points. We need to solve the equation,
$$\left[\frac{(\phi_{i,j}-\phi_{i,j}^{xmin})}{\Delta x}\right]^2+\left[\frac{(\phi_{i,j}-\phi_{i,j}^{ymin})}{\Delta y}\right]^2=\frac{1}{F_{i,j}^2},\eqno(5) $$
For simplicity of explanation, we assume a square mesh, so $\Delta x = \Delta y = \delta$. The local solver algorithm for 2D is,\\

\small
\begin{center}
\begin{minipage}{0.9\textwidth}
\noindent\rule[0.15\baselineskip]{\textwidth}{1.5pt}
$\textbf{Algorithm 1.1}$: $2D$ $update$ $function$ $at$ $\phi_{i,j}$ \\
\noindent\rule[0.15\baselineskip]{\textwidth}{1pt}
$\textbf{Find the neighboring points it depends on:}$\\
\hspace*{2em}$\phi^\text{minx} = \min(\phi_{i-1,j},\phi_{i+1,j})$\\
\hspace*{2em}$\phi^\text{miny} = \min(\phi_{i,j-1},\phi_{i,j-1})$\\
$\textbf{Update:}$\\
\hspace*{2em}$\textbf{if}$ $|\phi^\text{minx} -\phi^\text{miny}|>\sqrt 2\delta/f_{i,j}$: $\phi_{i,j}=\min(\phi^\text{minx} ,\phi^\text{miny})+\delta/f_{i,j}$\\
\hspace*{2em}$\textbf{else}$: $\phi_{i,j}=\frac{1}{2}\left [\phi^\text{minx} +\phi^\text{miny}+\sqrt{2\delta^2/f_{i,j}^2-(\phi^\text{minx} -\phi^\text{miny})^2} \quad\right]$\\
\noindent\rule[0.15\baselineskip]{\textwidth}{1.5pt}\\\\
\end{minipage}
\end{center}
\normalsize

\noindent They are a result of solving Equation (5), which is a quadratic equation in terms of the unknown $\phi_{i,j}$. For the case where $\Delta x\neq\Delta y$ Algorithm 1.1 is,\\

\small
\begin{center}
\begin{minipage}{0.9\textwidth}
\noindent\rule[0.15\baselineskip]{\textwidth}{1.5pt}
$\textbf{Algorithm 1.2}$: $2D$ $update$ $function$ $for$ $\phi_{i,j}$ \\
\noindent\rule[0.15\baselineskip]{\textwidth}{1pt}
$\textbf{Find the neighbour points it depends on:}$\\
\hspace*{1em}$\phi^\text{minx} = \min(\phi_{i-1,j},\phi_{i+1,j})$\\
\hspace*{1em}$\phi^\text{miny} = \min(\phi_{i,j-1},\phi_{i,j-1})$\\
$\textbf{Update:}$\\
\hspace*{1em}$\textbf{if}$  $(\phi^\text{minx}-\phi^\text{miny})>\sqrt{\frac{\Delta x^2 + \Delta y^2}{f_{i,j}^2}}$: $\phi_{i,j}=\phi^\text{miny}+\Delta y/f_{i,j}$\\
\hspace*{1em}$\textbf{else if}$  $(\phi^\text{minx}-\phi^\text{miny})<-\sqrt{\frac{\Delta x^2 + \Delta y^2}{f_{i,j}^2}}$: $\phi_{i,j}=\phi^\text{minx}+\Delta x/f_{i,j}$\\
\hspace*{1em}$\textbf{else}$: \\
\hspace*{1em}$\phi_{i,j}=\frac{1}{\Delta x^2 + \Delta y^2}\left[\phi^\text{minx}\Delta y^2+\phi^\text{miny}\Delta x^2 + \Delta x\Delta y\sqrt{\frac{\Delta x^2 + \Delta y^2}{f_{i,j}^2}-(\phi^\text{minx}-\phi^\text{miny})^2}\quad\right]$\\
\noindent\rule[0.15\baselineskip]{\textwidth}{1.5pt}\\\\
\end{minipage}
\end{center}
\normalsize

\noindent In 3D, the equation we need to solve for $\phi_{i,j,k}$ is,\\\\
$$\left[\frac{(\phi_{i,j,k}-\phi_{i,j,k}^{xmin})^+}{\Delta x}\right]^2+\left[\frac{(\phi_{i,j,k}-\phi_{i,j,k}^{ymin})^+}{\Delta y}\right]^2+\left[\frac{(\phi_{i,j,k}-\phi_{i,j,k}^{zmin})^+}{\Delta z}\right]^2=\frac{1}{F_{i,j,k}^2}.\eqno(6) $$\\\\\\\\\\
For simplicity of explanation, we assume a square mesh, so $\Delta x = \Delta y = \Delta z = \delta$. The local solver algorithm for 3D is,\\

\small
\begin{center}
\begin{minipage}{0.9\textwidth}
\noindent\rule[0.15\baselineskip]{\textwidth}{1.5pt}
$\textbf{Algorithm 1.3}$: $3D$ $update$ $function$ $at$ $\phi_{i,j,k}$ \\
\noindent\rule[0.15\baselineskip]{\textwidth}{1pt}
$\textbf{Find the neighboring points it depends on:}$\\
\hspace*{2em}$\phi^\text{minx} = \min(\phi_{i-1,j,k},\phi_{i+1,j,k})$\\
\hspace*{2em}$\phi^\text{miny} = \min(\phi_{i,j-1,k},\phi_{i,j-1,k})$\\
\hspace*{2em}$\phi^\text{minz} = \min(\phi_{i,j,k-1},\phi_{i,j,k-1})$\\
\hspace*{2em}$a_1$, $a_2$, $a_3$ $\longleftarrow$ sort$(\phi^\text{minx},\phi^\text{miny} ,\phi^\text{minz})$, where $a_1\leq a_2\leq a_3$.\\
$\textbf{Update:}$\\
\hspace*{2em}$\textbf{if}$ $|a_1-a_3|<\delta$: \\
\hspace*{4em}$\phi_{i,j,k}=\frac{1}{6}\left[\quad2\sum_{i=1}^3 a_i+\sqrt{4(\sum_{i=1}^3 a_i)^2-12(\sum_{i=1}^3a_i^2-\delta^2/f_{i,j,k}^2)}\quad\right] $\\
\hspace*{2em}$\textbf{else if}$ $|a_1-a_2|<\delta$: $\phi_{i,j,k}=\frac{1}{2}\left[\quad a_1+a_2+\sqrt{2\delta^2/f_{i,j,k}^2-(a_1-a_2)^2}\quad \right]$\\
\hspace*{2em}$\textbf{else}$: $\phi_{i,j,k}=a_1+\delta/f_{i,j,k}$\\
\noindent\rule[0.15\baselineskip]{\textwidth}{1.5pt}\\
\end{minipage}
\end{center}
\normalsize

\subsubsection{Main Algorithm}
There are three steps in the main algorithm. The first step is initialization where we set the values of grid points determined by the boundary conditions of the given problem and set the remaining values to $+\infty$. Then, we add the neighboring points to the boundary into the $Active$ set. The second step is to update the points in the $Active$ set until the set is empty. The last step is the remedy step. In this step, we first recalculate all the points again and add the points whose values change to the $Remedy$ set. Then, we update all points in the $Remedy$ set as well as their neighbors. After that, any point whose value changes is added to the $Remedy$ set and any point whose value does not change is removed from the Remedy set. The $Remedy$ step will end when the $Remedy$ set is empty.\\\\
The 2D version of the algorithm is shown below, the 3D version is analogous.
\small
\begin{center}
\begin{minipage}{0.9\textwidth}
\noindent\rule[0.15\baselineskip]{\textwidth}{1.5pt}
$\textbf{Algorithm 2}$: $Improved$ $Fast$ $Iterative$ $method$ \\
\noindent\rule[0.15\baselineskip]{\textwidth}{1pt}
$\textbf{STEP 1 : Variable Define \& Initialization}$\\
$\bullet$ Define a set $Source$ including all points assigned an initial value by the problem to be solved.\\
$\bullet$ Assign all values $\phi_{i,j}$ in the Source set. All other points are assigned the value $+\infty$.\\
$\bullet$ Define the $Active$ set to be the set of points adjacent to the points in the $Source$ set. \\\\
$\textbf{STEP 2 : Update}$\\
\textbf{while} the Active set is not empty:\\
$\bullet$ Update the point values in the $Active$ set:\\
\hspace*{2em} \textbf{if} $\mathbf{x}_{i,j}\in Active$: Calculate $\tilde{\phi}_{i,j}$ based on \textbf{Algorithm 1}\\
\hspace*{4em} \textbf{if} $\tilde{\phi}_{i,j} = \phi_{i,j}$:\\
\hspace*{4em} Remove $\mathbf{x}_{i,j}$ from the $Active$ set and add it to the $Source$ set.\\
\hspace*{6em} \textbf{if} any neighbor point $\mathbf{x}_{i\pm1,j}, \mathbf{x}_{i,j\pm1}\notin Active\cup Source$\\
\hspace*{6em} Add it to the $Active$ set.\\
\hspace*{2em} \textbf{else}: $\phi_{i,j}=\tilde{\phi}_{i,j}$.\\\\
$\textbf{STEP 3: Remedy}$\\
$\bullet$ Define the $Remedy$ set to be the set of points that don't satisfy Equation (5):\\
\hspace*{2em} Calculate $\tilde{\phi}_{i,j}$ for all grid points.\\
\hspace*{4em} \textbf{if} $\tilde{\phi}_{i,j} \neq \phi_{i,j}$: then add $\mathbf{x}_{i,j}$ to the Remedy set.\\
\textbf{do}:\\
\hspace*{2em} \textbf{if} the point $\mathbf{x}_{i,j}\in Remedy$: \\
\hspace*{4em} Calculate $\tilde{\phi}_{i,j}$ based on \textbf{Algorithm 1} \\
\hspace*{4em} \textbf{if} $\tilde{\phi}_{i,j} < \phi_{i,j}$ :\\
\hspace*{6em} $\phi_{i,j} = \tilde{\phi}_{i,j}$.\\
\hspace*{6em} \textbf{if} any neighbor $\mathbf{x}_{i\pm1,j}, \mathbf{x}_{i,j\pm1}\notin Remedy$: \\
\hspace*{6em} Add it to the $Remedy$ set.\\
\hspace*{4em} \textbf{else}: \\
\hspace*{6em} Remove the point $\mathbf{x}_{i,j}$ from the $Remedy$ set.\\
\textbf{while} the Remedy set is not empty:\\
\noindent\rule[0.15\baselineskip]{\textwidth}{1.5pt}
\end{minipage}
\end{center}
\normalsize

\subsection{Parallel Implementation}
\subsubsection{Parallel for $\phi$ Update Step}
Each iteration, we calculate the $\phi$ values of the points in the $Active$ set in parallel. On a GPU, each thread calculates $\phi$ at one point each iteration. After all threads finish computing $\phi$ in a given iteration, threads start updating the $Active$ and $Source$ set for the grid point each thread updated. That completes one iteration of the method.

\subsubsection{Parallel for Remedy Step}
Parallel computing for the remedy step is handled in the same way as the update step in Section 3.2.1.

\subsection{Compared with FIM}
Won-Ki Jeong and Ross Whitaker [6] proposed the fast iterative method. One property of the FIM is that most points require only a single update to converge, so it is inefficient to check the convergence of the neighboring points of the $Active$ list in every iteration. Even in parallel, each thread needs to update the point value of a point in the $Active$ list and then check the convergence of its neighbors. Threads cannot easily do these two jobs simultaneously.\\\\
The improved FIM increases the efficiency of the parallel implementation in two ways. The \textbf{Update Step} does not check for convergence of the neighbor points of the $Active$ set saving much time. Since most points require only a single update to converge due to the property of the FIM, only a few points will have errors after the update step. The \textbf{Remedy Step} can be done so quickly in parallel because generally relatively few points need to be remedied.\\\\\\

\section{Results}

\subsection{Examples}
For the numerical experiments, in order to prove the efficiency of the parallel version of the algorithm, we implemented the serial version of the fast marching method and the fast sweeping method on a single core Intel(R) Xeon(R) CPU. We implemented the parallel version of the fast iterative method and the improved fast iterative method using an NVIDIA Tesla P100 GPU. Some examples we used are the same examples from paper [6] for purposes of comparison.

\subsubsection*{Example 1}
This example shows a simple circle initial boundary $\Gamma$ expanding with a constant speed function $F=1$. Note that when $F=1$, there is no need to maintain a heap queue for the fast marching method. Fig 1 shows the isocontours of the computed solution for $\phi$. Since the speed function is a constant equal to $1$, the isocontouor with value $\phi=t$ is a circle with radius $t+R_{initial}$. \\
$$\includegraphics[scale=0.2]{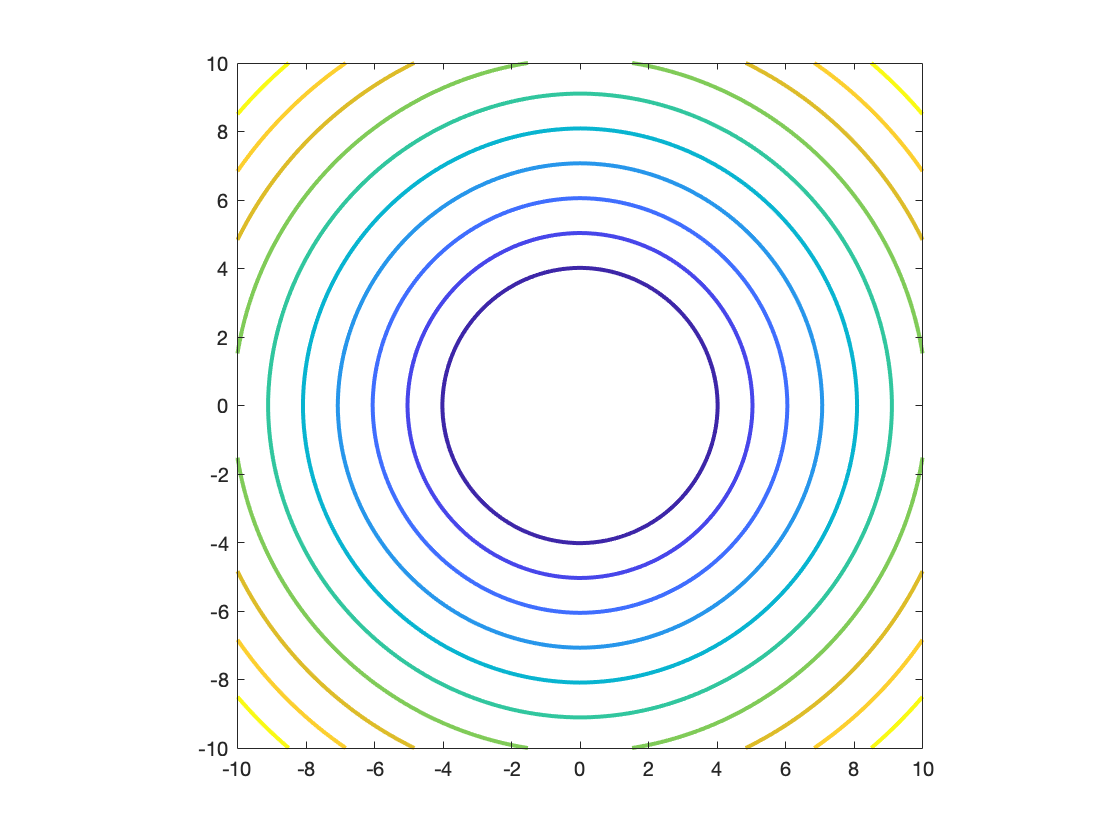}$$
\small
\begin{center}
\begin{minipage}{0.7\textwidth}
Figure 1. Plot of the solution for the case of $F=1$ in Equation (1) where the initial condition is a circle of radius $X=3$. Note that all methods produce exactly the same plots, so only one plot is shown. (Since these four algorithms get exactly the same solution, we only show one figure here.)\\\\
\end{minipage}
\end{center}
\normalsize
\begin{center}
\begin{tabular}{cccccc}
\\
\toprule
N&$200^2$&$500^2$&$1000^2$&$2000^2$&$4000^2$\\
\midrule
FMM&0.00780&0.13783&0.2918&1.22073&3.7991\\
FSM&0.00610&0.16702&0.36897&1.49810&4.84172\\
FIM&0.00164&0.00428&0.08834&0.32604&1.67655\\
\textbf{iFIM}&\textbf{0.00141}&\textbf{0.00384}&\textbf{0.03721}&\textbf{0.11642}&\textbf{0.80981}\\
\bottomrule
\end{tabular}
\end{center}
\small
\begin{center}
\begin{minipage}{0.7\textwidth}
Table 1. Elapsed Time vs. Number of Grid Points. Timing data of the solution for the case of $F=1$ in Equation (1) where the initial conditions were a circle of radius $X=3$. The table shows the elapsed time in seconds of the four algorithms against the number of gridpoints $N$.
\end{minipage}
\end{center}
\normalsize
\begin{center}
$\includegraphics[scale = 0.6]{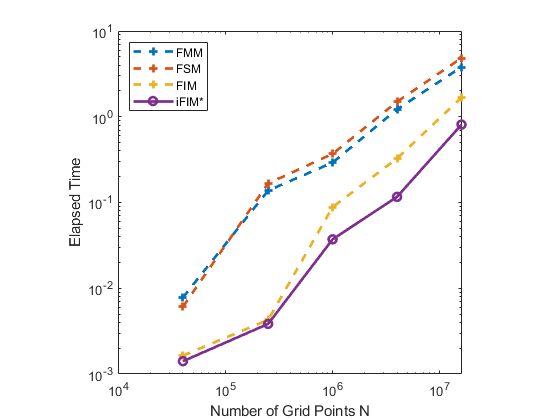}$
\small
\begin{minipage}{0.7\textwidth}
Figure 2. Elapsed Time vs. Number of Grid Points. Plot of the data of the solution for the case of $F=1$ in Equation (1) where the initial conditions were a circle of radius $X=3$. The graph shows the elapsed times in seconds of the four algorithms against the number of gridpoints $N$.\\\\
\end{minipage}
\end{center}
\normalsize
\noindent Table 1 and Fig 2 show the elapsed time in seconds for each algorithm. We see that the improved fast iterative method is consistently faster than the other methods. It is because there is no need for the remedy step for this simple initial boundary.\\\\\\

\subsubsection*{Example 2}

This example shows an initial boundary $\Gamma$ that is the combination of two circles expanding with a constant speed function $F = 1$. For this example, there is still no need to maintain a heap queue in the FMM. Fig 3 shows the solution of $\phi$. The initial condition is two circles of radius $X_1 = 3$ and $X_2=1.5$, centered at $(2,-5)$ and $(-2,5)$. Table 2 and Fig 4 show the elapsed time in seconds for each algorithm.
\begin{center}
$$\includegraphics[scale=0.2]{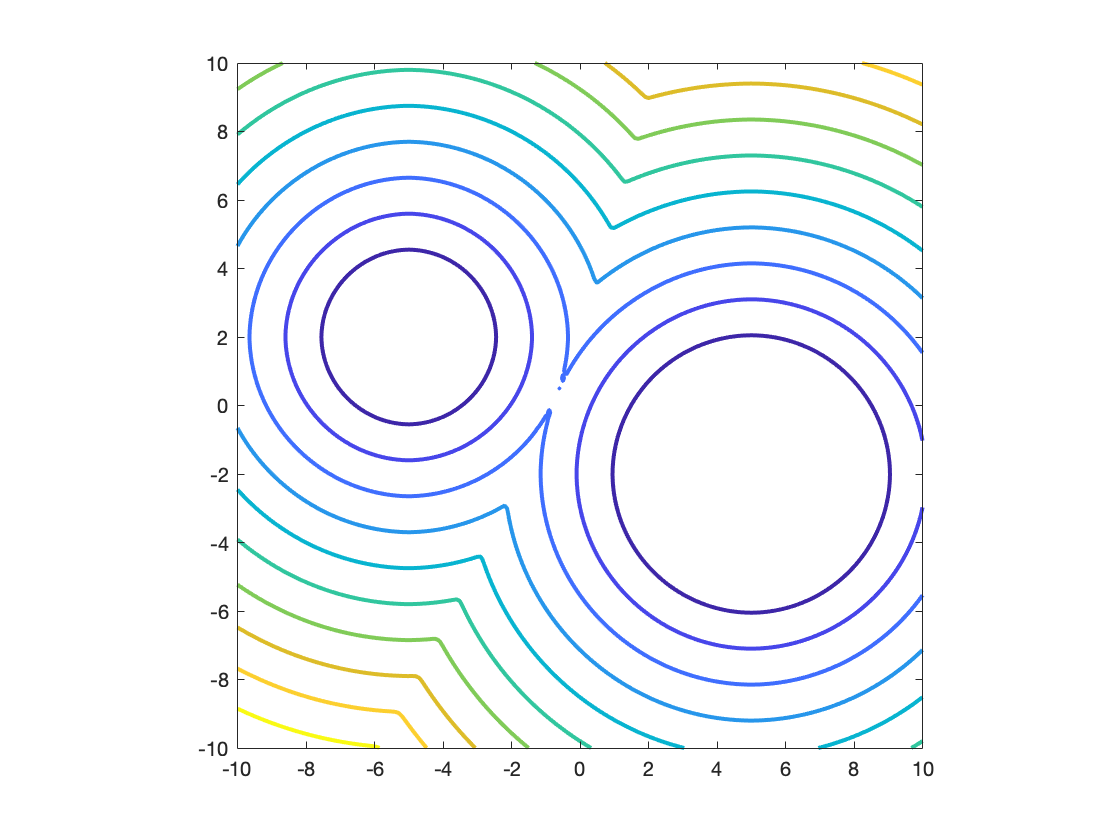}$$
\small
\begin{minipage}{0.7\textwidth}
Figure 3. Plot of the solution for the case of $F=1$ in Equation (1) where the initial condition is two circles of radius $X_1=3$ and $X_2=1.5$, centered at $(2,-5)$ and $(-2,5)$ respectively. Note that all methods produce exactly the same plots, so only one plot is shown. (Since these four algorithms get exactly the same solution, we only show one figure here.)\\\\\\\\\\
\end{minipage}
\end{center}
\normalsize

\begin{center}
\begin{tabular}{l|cccccc}
\toprule
N&$200^2$&$500^2$&$1000^2$&$2000^2$&$4000^2$\\
\midrule
FMM&0.02740&0.15031&0.33081&1.10910&3.5799\\
FSM&0.03551&0.16731&0.42476& 1.30234&4.75479\\
FIM&0.00118&0.00319&0.08794&0.63421&1.86989\\
\textbf{iFIM}&\textbf{0.00097}&\textbf{0.00329}&\textbf{0.0539}&\textbf{0.35359}&\textbf{1.09434}\\
\bottomrule
\end{tabular}
\end{center}
\begin{center}
\small
\begin{minipage}{0.7\textwidth}
Table 2. Elapsed Time vs. Number of Grid Points. Timing data of the solution for the case of $F=1$ in Equation (1) where the initial condition is two circles of radius $X_1=3$ and $X_2=1.5$, centered at $(2,-5)$ and $(-2,5)$ respectively. The table shows the elapsed time in seconds of the four algorithms against the number of gridpoints $N$.\\
\end{minipage}
\end{center}
\normalsize

\begin{center}
$$\includegraphics[scale = 0.6]{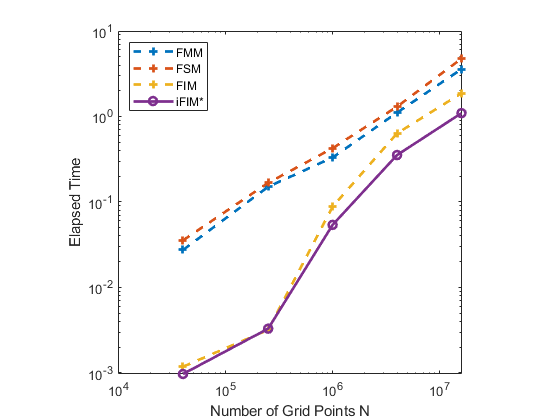}$$
\small
\begin{minipage}{0.7\textwidth}
Figure 4. Elapsed Time vs. Number of Grid Points. Plot of the data of the solution for the case of $F=1$ in Equation (1) where the initial condition is two circles of radius $X_1=3$ and $X_2=1.5$, centered at $(2,-5)$ and $(-2,5)$ respectively. The graph shows the elapsed time in seconds of the four algorithms against the number of gridpoints $N$.\\\\
\end{minipage}
\end{center}
\normalsize

\noindent We see that both parallel methods required more time compared with that of Example~1 because the remedy step is required for this example.\\\\\\

\subsubsection*{Example 3}
This example shows an initial boundary $\Gamma$ that is two small circles of radius $X_1=X_2=0.5$, centered at $(-5,5)$ and $(5,-5)$, expanding with a variable speed function $F$. For this situation, it is necessary to maintain a heap queue in the FMM. Fig 5 shows the level sets of the solution $\phi$ where the speed function is,
$$F = \left\{\begin{array}{lll}
1&  & (x,y) \in[3,7]\times[3,7] ,\\
& \\
0.01&  & \text{otherwise}  
 \end{array}\right.\eqno(7)$$
We implement it on a grid where $x$ and $y$ range from $-10$ to $10$. Table 3 and Fig 6 show the elapsed time in seconds for each algorithm.

\clearpage
\begin{center}
$$\includegraphics[scale=0.15]{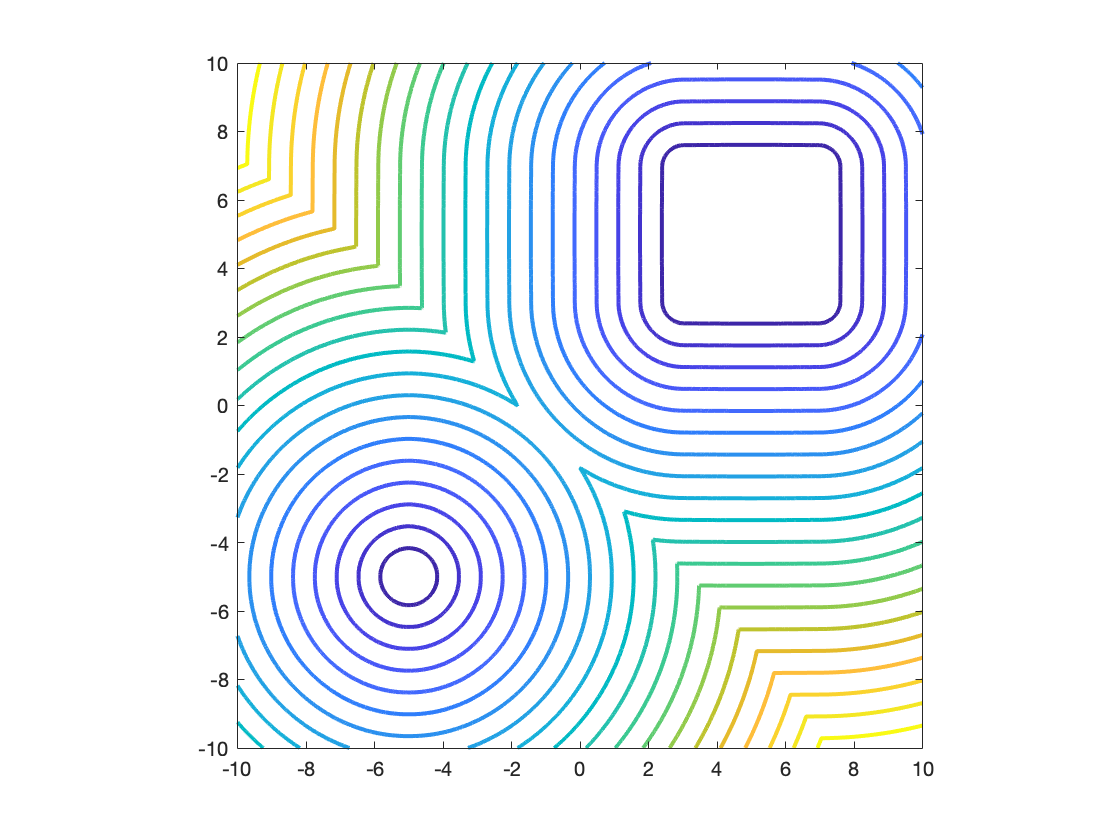}$$
\small
\begin{minipage}{0.7\textwidth}
Figiure 5. Plot of the solution for the case of $F$ in Equation (7) where the initial condition is two circles of radius $X_1=X_2=0.5$, centered at $(-5,5)$ and $(5,-5)$ respectively. Note that all methods produce exactly the same plots, so only one plot is shown. (Since these four algorithms get exactly the same solution, we only show one figure here.)\\\\
\end{minipage}
\end{center}
\normalsize

\begin{center}
\begin{tabular}{l|cccccc}
\toprule
N&$200^2$&$500^2$&$1000^2$&$2000^2$&$4000^2$\\
\midrule
FMM&0.04512&0.39215&1.39354&4.98564&20.9888\\
FSM&0.03959&0.17895&0.43589&1.35517&5.12889\\
FIM&0.00125&0.00410&0.08954&0.78251&2.05581\\
\textbf{iFIM}&\textbf{0.00089}&\textbf{0.00341}&\textbf{0.0499}&\textbf{0.36897}&\textbf{1.48971}\\
\bottomrule
\end{tabular}
\end{center}
\begin{center}
\small
\begin{minipage}{0.7\textwidth}
Table 3. Elapsed Time vs. Number of Grid Points. Timing data of the solution for the case of $F$ in Equation (7) where the initial condition is two circles of radius $X_1=X_2=0.5$, centered at $(-5,5)$ and $(5,-5)$ respectively. The table shows the elapsed time in seconds of the four algorithms against the number of gridpints $N$.\\\\
\end{minipage}
\end{center}
\normalsize

\begin{center}
$$\includegraphics[scale=0.5]{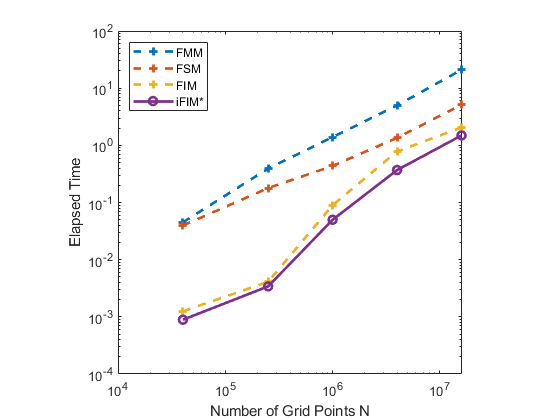}$$
\small
\begin{minipage}{0.7\textwidth}
Figure 6. Elapsed Time vs. Number of Grid Points. Plot of the data of the solution for the case of $F$ in Equation (7) where the initial condition is two circles of radius $X_1=X_2=0.5$, centered at $(-5,5)$ and $(5,-5)$. The graph shows the elapsed time in seconds of the four algorithms against the number of the gridpoints $N$.\\\\
\end{minipage}
\end{center}
\normalsize

\noindent We see that the elapsed time for FMM is much longer when a heap queue is required.\\\\

\subsubsection*{Example 4}
For a more realistic example, we compute the shortest path between two points on a map. We set the speed function $F=0$ on the grid points which are in the barrier area while setting $F=1$ everywhere else,
$$F = \left\{\begin{array}{lll}
0&  & (x,y) \in \text{MAP}_{barrier} \\
& \\
1&  & \text{otherwise}  
 \end{array}\right.\eqno(8)$$
 and set the starting point as the initial boundary $\Gamma$ located at (490, 245). Fig 7 shows the isocontours of the soluton $\phi$.
 
\begin{center}
$$\includegraphics[scale=0.22]{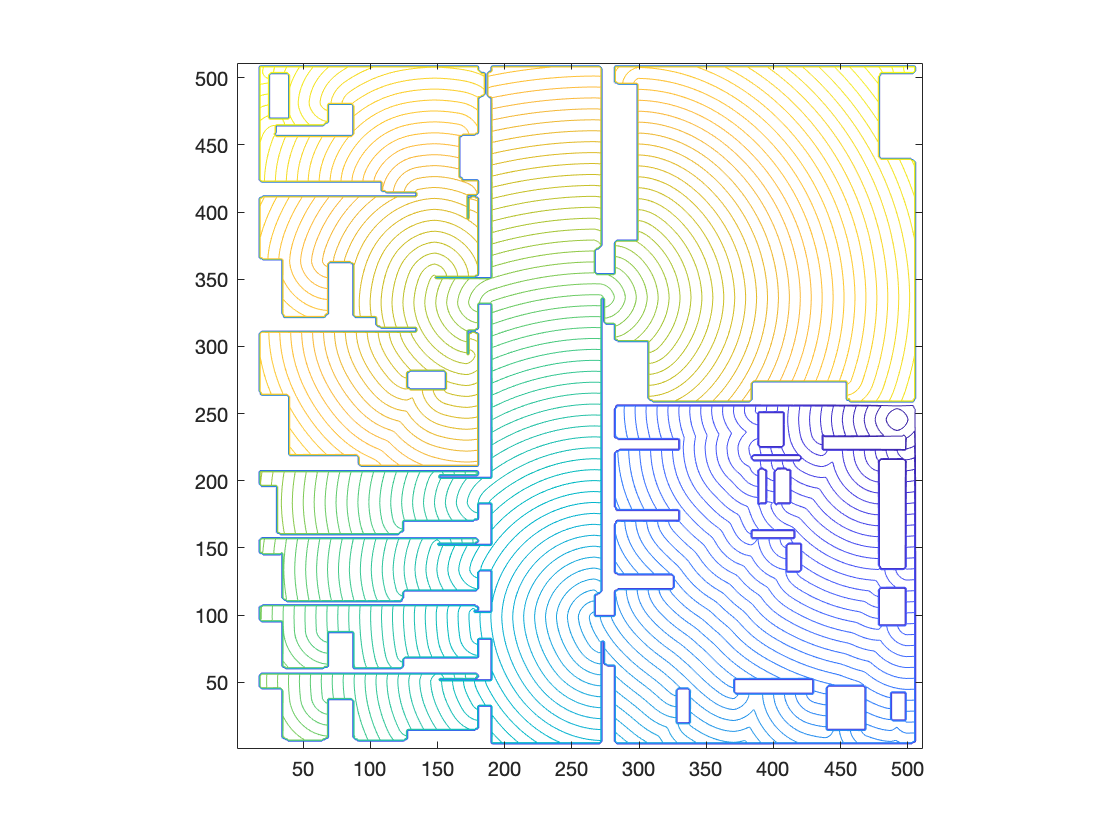}$$
\small
\begin{minipage}{0.7\textwidth}
Figure 7. $\phi$ Map. Plot of the solution for the case of $F$ in Equation (8) where the initial condition is a point located at (490, 245). Note that all methods produce exactly the same plots, so only one plot is shown. (Since these four algorithms get exactly the same solution, we only show one figure here.)\\\\
\end{minipage}
\end{center}
\normalsize

\noindent Then, we use the gradient descent algorithm on the $\phi$ map to find the shortest path between a destination point located at (50, 475) and the starting point [13]. Fig 8 shows the shortest path,

\begin{center}
$$\includegraphics[scale=0.22]{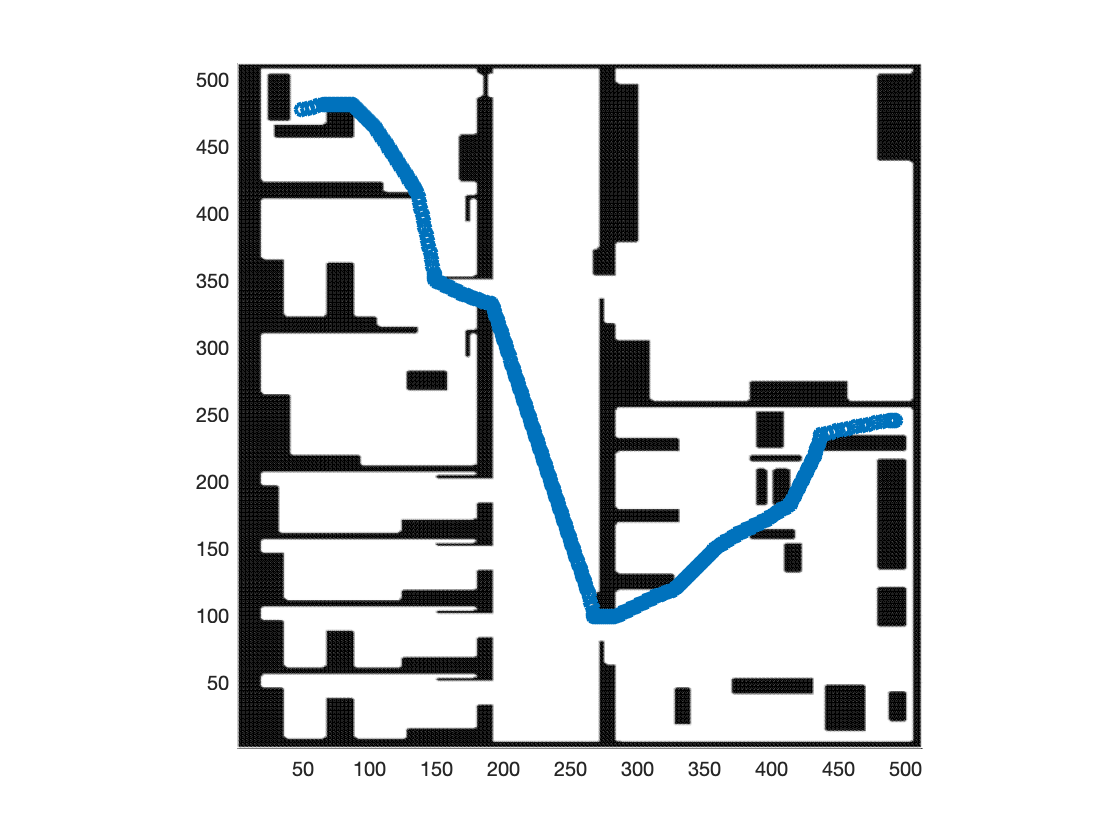}$$
\small
\begin{minipage}{0.7\textwidth}
Figure 8. Shortest Path. Shortest path between the point (490, 245) and (50, 475) by using the gradient descent algorithm on the $\phi$ map.\\\\
\end{minipage}
\end{center}
\normalsize

\begin{center}
\begin{tabular}{lccccc}
\toprule
N&$256^2$&$512^2$&$1024^2$&$2048^2$\\
\midrule
FMM&0.03856&0.13094&0.437167&2.01458\\
FSM&0.14650&0.37300&1.37545&6.56971\\
FIM&0.00640&0.02234&0.32843&1.42556\\
\textbf{iFIM}&\textbf{0.00484}&\textbf{0.01699}&\textbf{0.26843}&\textbf{0.98356}\\
\bottomrule
\end{tabular}
\end{center}
\begin{center}
\small
\begin{minipage}{0.7\textwidth}
Table 4. Elapsed Time vs. Number of Grid Points. Timing data of the solution for the case of $F$ in Equation (8) where the initial condition is a point located at (490, 245). The table shows the elapsed time in seconds of the four algorithms against the number of gridpoints $N$.
\end{minipage}
\end{center}
\normalsize

\begin{center}
$$\includegraphics[scale=0.6]{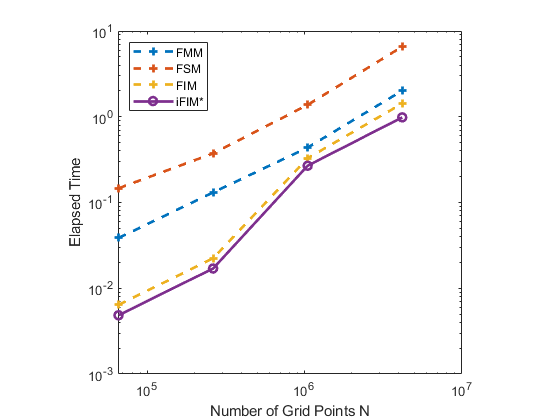}$$
\small
\begin{minipage}{0.7\textwidth}
Figure 9. Elapsed Time vs. Number of Grid Points. Plot of the timing data of the solution for the case of $F$ in Equation (8) where the initial condition is a point located at (490, 245). The graph shows the elapsed time in seconds of the four algorithms against the number of gridpoints $N$.\\\\
\end{minipage}
\end{center}
\normalsize

\noindent Table 4 and Fig 9 shows the elapsed time in seconds for each algorithm. We see from the results that the elapsed time for the FSM is much longer. It is because it cannot get all the grid points to converge in only four sweeps due to the barriers on the map.\\\\\\

\subsubsection*{Example 5}
This example shows an initial boundary $\Gamma$ of a small circle expanding with a variable speed function $F$. For this situation, it is necessary to maintain a heap queue in FMM due to the spatially varying $F$ function. Fig 10 shows the level sets of $\phi$ for the speed function,
$$F = \sin^2(x)+\sin^2(y)+\epsilon,\quad \epsilon = 0.3,\eqno(9)$$
and we implement it on a grid where $x$ and $y$ range from $-10$ to $10$. Table 5 and Fig 11 show the elapsed time in seconds for each algorithm.

\begin{center}
$$\includegraphics[scale=0.2]{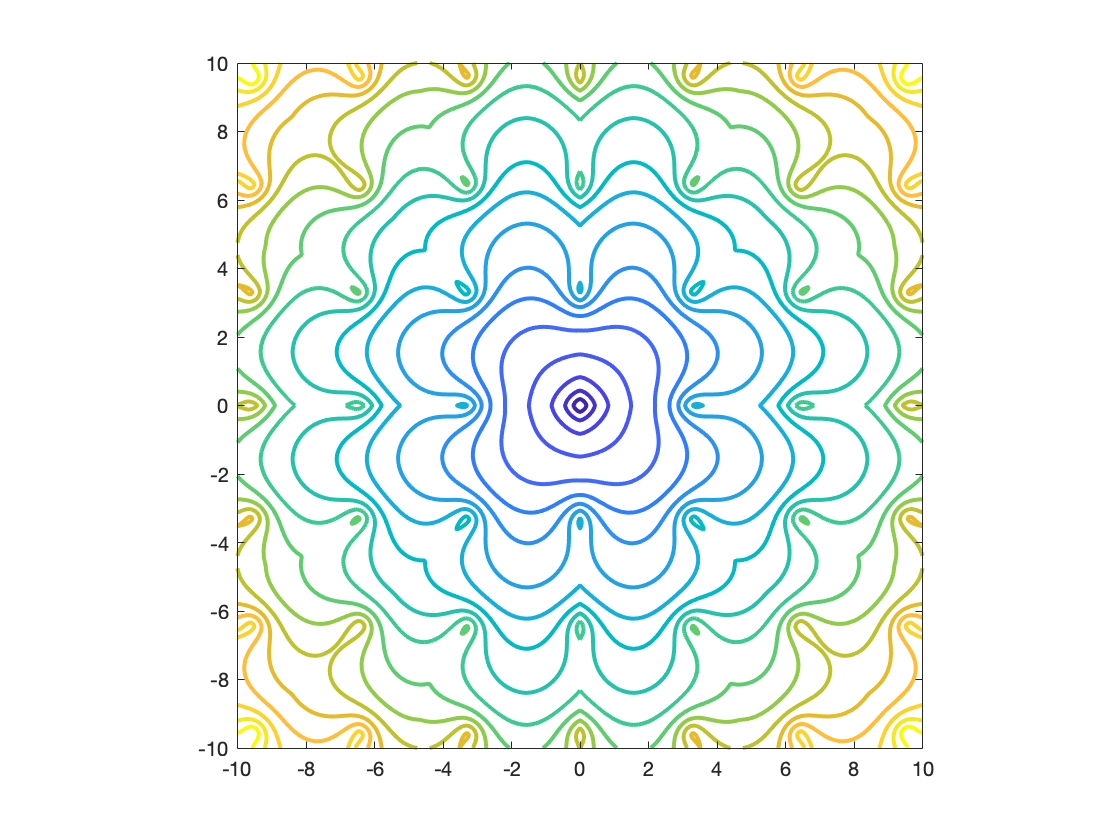}$$
\small
\begin{minipage}{0.7\textwidth}
Figure 10. Plot of the solution for the case of $F$ in Equation~(9) where the initial condition is a circle of radius $X=0.5$ centered at the origin. Note that all methods produce exactly the same plots, so only one plot is shown. (Since these four algorithms get exactly the same solution, we only show one figure here.)\\\\
\end{minipage}
\end{center}
\normalsize

\begin{center}
\begin{tabular}{l|cccccc}
\toprule
N&$200^2$&$500^2$&$1000^2$&$2000^2$&$4000^2$\\
\midrule
FMM&0.05961&0.49815&1.41358&5.28564&22.9005\\
FSM&0.04259&0.35154&0.48931&1.52147&5.92456\\
FIM&0.00125&0.00410&0.08954&0.78251&2.45581\\
\textbf{iFIM}&\textbf{0.00119}&\textbf{0.00341}&\textbf{0.0689}&\textbf{0.56897}&\textbf{1.098971}\\
\bottomrule
\end{tabular}
\end{center}
\begin{center}
\small
\begin{minipage}{0.7\textwidth}
Table 5. Elapsed Time vs. Number of Grid Points. Timing data of the solution for the case of $F$ in Equation~(9) where the initial condition is a circle of radius $X=0.5$ centered at the origin. The table shows the elapsed time in seconds of the four algorithms against the number of $N$.\\\\
\end{minipage}
\end{center}
\normalsize

\begin{center}
$$\includegraphics[scale=0.5]{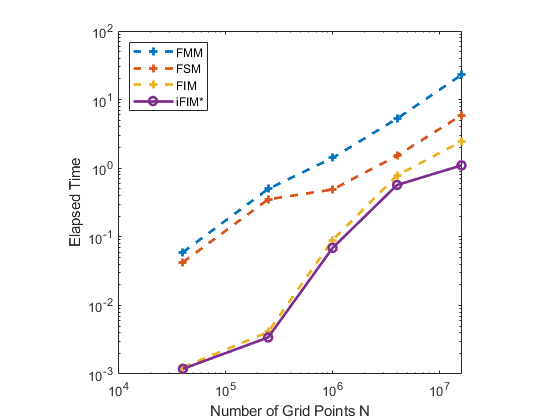}$$
\small
\begin{minipage}{0.7\textwidth}
Figure 11. Elapsed Time vs. Number of Grid Points. Plot of the timing data of the solution for the case of $F$ in Equation~(9) where the initial condition is a circle of radius $X=0.5$ centered at the origin. The graph shows the elapsed time in seconds of the four algorithms against the number of gridpoints $N$.\\\\
\end{minipage}
\end{center}
\normalsize

\noindent We see from the examples that when $F$ is variable the elapsed time for FMM is always much longer because a heap sort queue is required.\\\\
\clearpage

\subsection{Analysis}
\subsubsection*{GPU Performance}
Though running time is affected by a lot of factors such as memory, communication in GPUs, single core performance, the parallel algorithms consistently run faster compared to the serial version. The most time consuming part of the numerical method for the Eikonal equation is solving the quadratic equations (\textbf{Algorithm 1}). In the serial algorithm, we implement the local solver point by point so the iteration number $I_n \propto N$ where $N$ denotes the total number of grid points. In the parallel algorithm, we implement the local solver on many points simultaneously. Thus, the implementation of the parallel version for the Eikonal equation can be much faster than that of the serial version of the algorithm.

\subsubsection*{Elapsed Time and Thread Number}
It should be noted that the results can depend significantly on the number of threads assigned to each block when invoking the kernels. Fig 12 shows the elapsed time for Example 2 as the number of threads per block varies for $N=2000^2$, We see that tuning the number of threads per block is important to maintain peak efficiency. \\

\begin{center}
$$\includegraphics[scale=0.5]{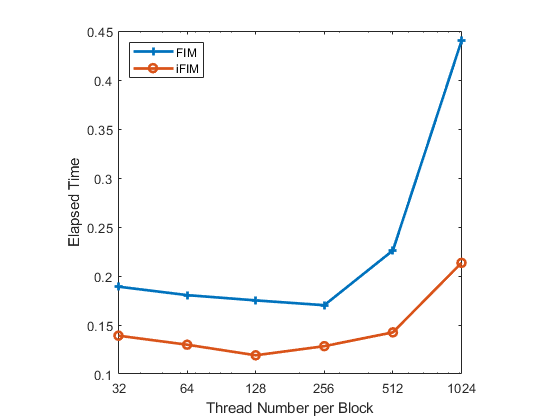}$$
\small
\begin{minipage}{0.7\textwidth}
Figure 12. Number of Threads per Block $vs.$ Elapsed Time of FIM and iFIM in Example 2.\\\\
\end{minipage}
\end{center}
\normalsize

\subsubsection*{Comparison between FIM and iFIM}
According to the results above, the improved fast iterative method consistently ran faster than the fast iterative method in parallel. This is because most of the grid points converge after only a single update as Won-Ki Jeong and Ross Whitaker [6] observed. For Example 1, we see that all the points are convergent after the Update Step so that no remedy is needed as shown shown in Fig 13,

\begin{center}
$$\includegraphics[scale=0.3]{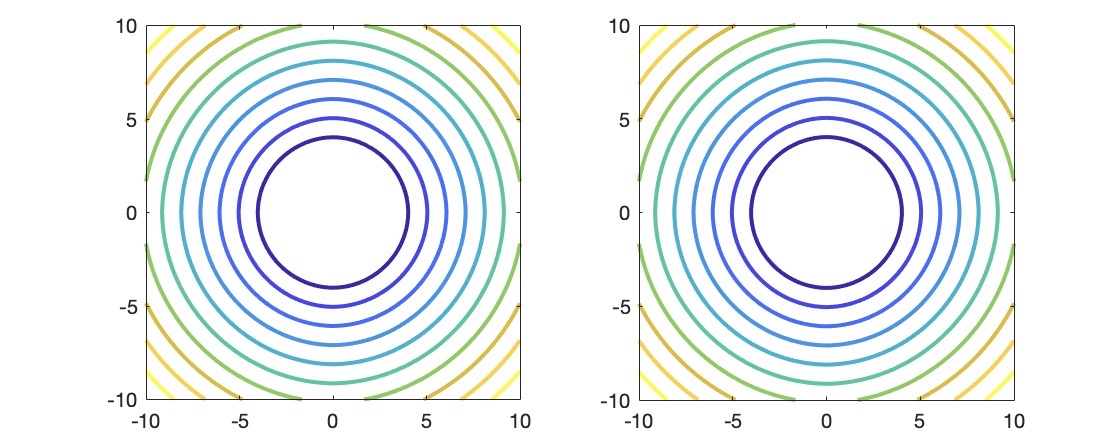}$$
\small
\begin{minipage}{0.7\textwidth}
Figure 13. Exact Solution (left) $vs.$ Solution after the \textbf{Update Step} (right)\\\\
\end{minipage}
\end{center}
\normalsize

\noindent For Example 2 and Example 5, some of the grid points require correction after the Update Step as shown in Fig 14 and Fig 15,

\begin{center}
$$\includegraphics[scale=0.3]{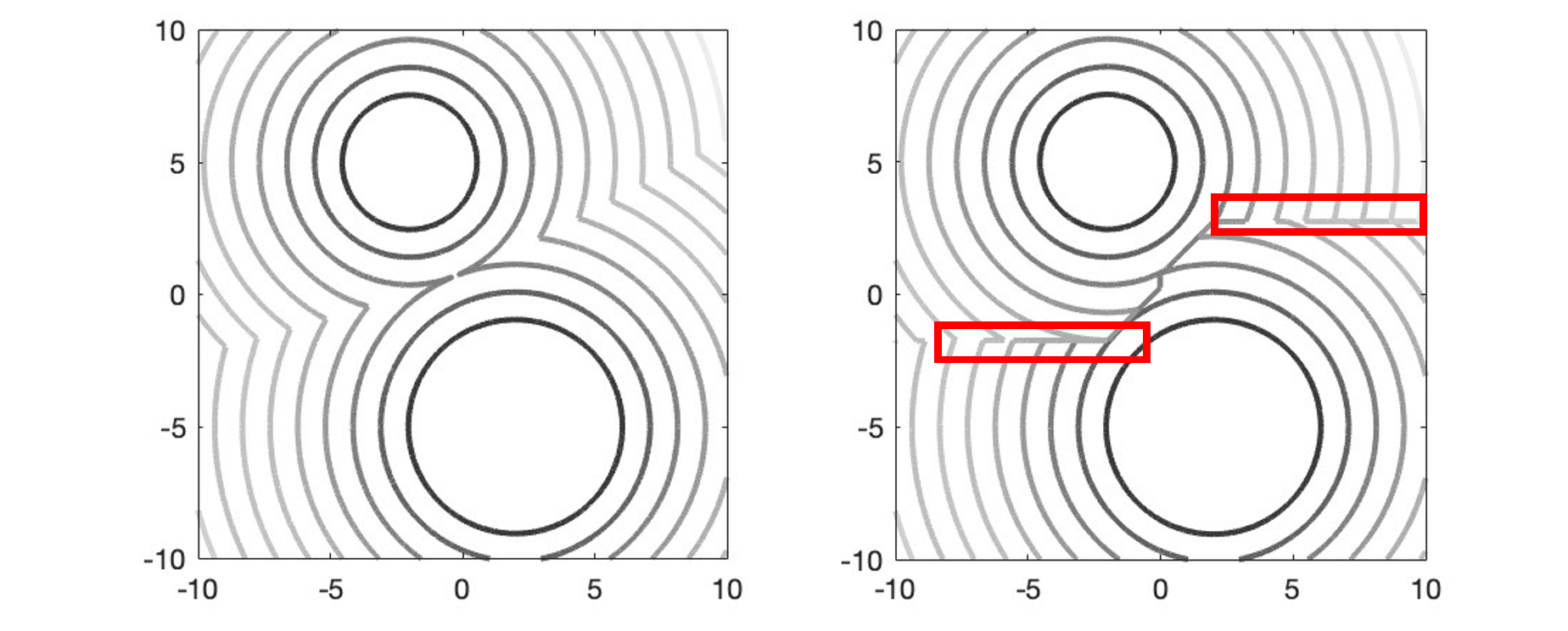}$$
\small
\begin{minipage}{0.7\textwidth}
Figure 14. Exact Solution (left) $vs.$ Solution after the Update Step (right). The left figure shows the exact solution to Example 2 and the right figure shows the solution to Example 2 given by the update step in iFIM. The red boxes in the right figure show the errors after the update step in iFIM.\\\\
\end{minipage}
\end{center}
\normalsize

\begin{center}
$$\includegraphics[scale=0.3]{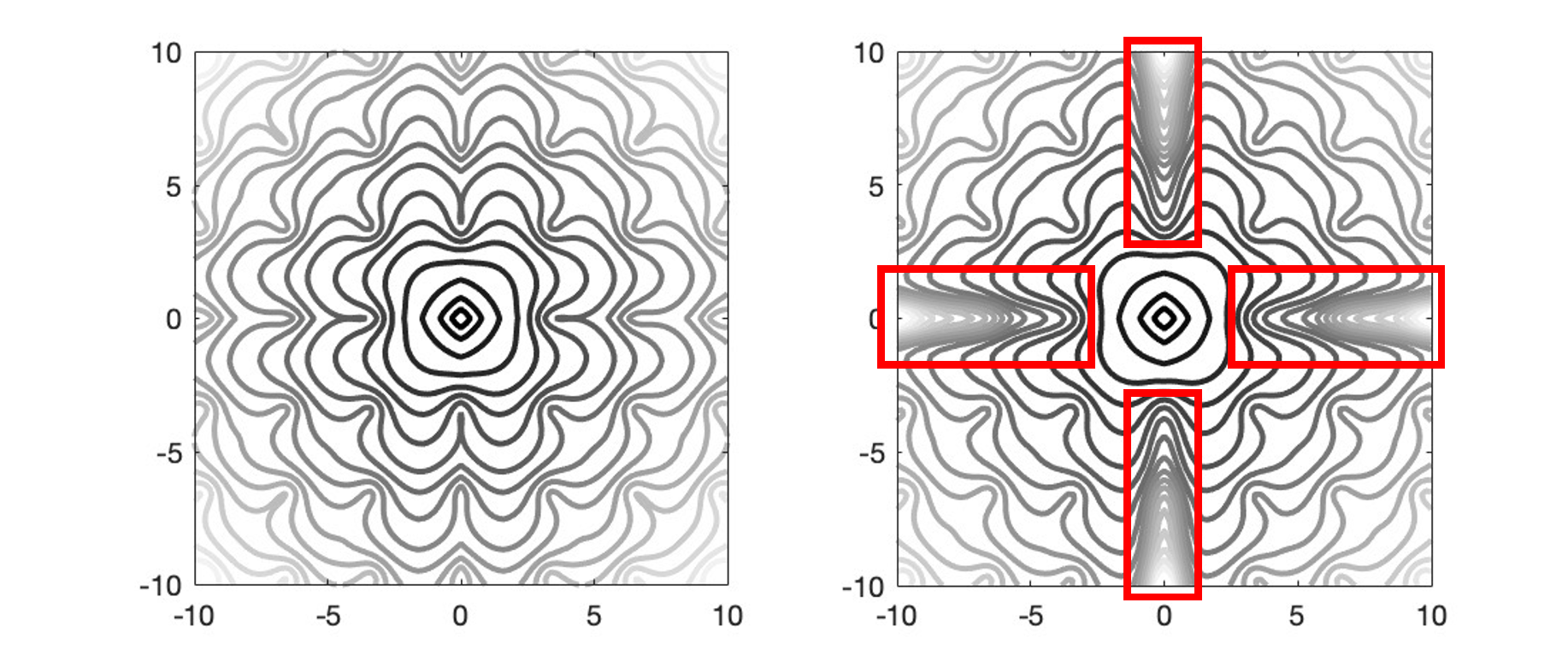}$$
\small
\begin{minipage}{0.7\textwidth}
Figure 15. Exact Solution (left) $vs.$ Solution after the Update Step(right). The left figure shows the exact solution to Example 5 and the right figure shows the solutions to Example 5 given by the update step in iFIM. The red boxes in the right figure show the errors after the update step in iFIM. Due to some errors in the right figure, the contour differences of these two figures are different.\\\\
\end{minipage}
\end{center}
\normalsize

\noindent In this case, some of the points have errors. In the parallel algorithm, we can quickly locate the error points and use the remedy step to get rid of the errors.\\\\
Therefore, based on the experimental results and the analysis, we see that the improved fast iterative method consistently runs faster than the fast iterative method.\\\\

\section{Summary and Future Work}
In this thesis, we propose an improved iterative algorithm for the Eikonal equation based on the fast iterative method. Compared with the fast iterative method, our method removed the convergence check step in each iteration of the update step and added the remedy step after update step. Our method improves the operation speed in GPUs according to our experiments and it is easy to implement. The method takes advantage of the ability of GPUs to increase the operation speed compared with the serial version of the algorithm. To support parallel computing, we keep two sets of points to be updated and remedy instead of maintaining a sorting data structure which is very expensive and does not support SIMD. \\\\
We developed this algorithm because it can take advantage of the GPUs. The remedy points checking step in our algorithm is a little expensive, especially when GPU computing resources are not used. Our future work is to improve the method further and to find a way to get rid of the remedy points checking step. More generally, developing a parallel method for Hamilton-Jacobi PDEs would also be another useful extension of this work.

\end{titlepage}

\clearpage
\begin{thebibliography}{9}
 \bibitem[1]{cite_key1} J. Sethian. (1999). Fast marching methods. SIAM Review, 41(2): 199 - 235.
 \bibitem[2]{cite_key2} J. Sethian. (2002). Level set methods and fast marching methods.  Cambridge
 University Press.
 \bibitem[3]{cite_key3} J. Sethian. (1996). A fast marching level set method for monotonically advancing fronts. In Proc. Natl. Acad. Sci., volume 93, pages 1591–1595.
 \bibitem[4]{cite_key4} H. Zhao. (2004). A fast sweeping method for eikonal equations. Mathematics of Computation, 74:603–627.
 \bibitem[5]{cite_key4} H. Zhao. (2007). Parallel implementations of the fast sweeping method. Journal of Computational Mathematics, 421-429.
 \bibitem[6]{cite_key4} Jeong, W. K., \& Whitaker, R. T. (2008). A fast iterative method for eikonal equations. SIAM Journal on Scientific Computing, 30(5), 2512-2534.
 \bibitem[7]{cite_key4} Fu, Z., Jeong, W. K., Pan, Y., Kirby, R. M., \& Whitaker, R. T. (2011). A fast iterative method for solving the eikonal equation on triangulated surfaces. SIAM Journal on Scientific Computing, 33(5), 2468-2488.
 \bibitem[8]{cite_key4} Rouy, E., \& Tourin, A. (1992). A viscosity solutions approach to shape-from-shading. SIAM Journal on Numerical Analysis, 29(3), 867-884.
 \bibitem[9]{cite_key4} Sethian, J. A., \& Popovici, A. M. (1999). 3-D traveltime computation using the fast marching method. Geophysics, 64(2), 516-523.
 \bibitem[10]{cite_key4} Rawlinson, N., \& Sambridge, M. (2004). Multiple reflection and transmission phases in complex layered media using a multistage fast marching method. Geophysics, 69(5), 1338-1350.
 \bibitem[11]{cite_key4} Rawlinson, N., \& Sambridge, M. (2005). The fast marching method: an effective tool for tomographic imaging and tracking multiple phases in complex layered media. Exploration Geophysics, 36(4), 341-350.
 \bibitem[12]{cite_key4} Garrido, S., Moreno, L., Abderrahim, M., \& Martin, F. (2006, October). Path planning for mobile robot navigation using voronoi diagram and fast marching. In 2006 IEEE/RSJ International Conference on Intelligent Robots and Systems (pp. 2376-2381). IEEE.
\end{thebibliography}
\end{document}